\documentclass[a4paper, reqno]{amsart}
\usepackage[latin1]{inputenc}
\usepackage{amsmath,amsthm,amssymb,amsfonts}
\usepackage{enumerate}

\theoremstyle{definition}
\newtheorem{definition}{Definition}[section]

\newtheorem{example}[definition]{Example}
\theoremstyle{remark}
\newtheorem{note}[definition]{Remark}

\theoremstyle{plain}
\newtheorem{lemma}[definition]{Lemma}

\newtheorem{prop}[definition]{Proposition}

\newtheorem{theorem}[definition]{Theorem}
\newtheorem{thm}[definition]{Theorem}
\newtheorem{corr}[definition]{Corollary}

\newtheorem{cor}[definition]{Corollary}

\newtheorem{question}[definition]{Question}

\newcommand{\Nat}{{\mathbb{N}}}

\newcommand{\C}{{\mathbb{C}}}
\newcommand{\R}{{\mathbb{R}}}

\def\setsuchas#1#2{\left\{\,{#1}\,\vrule\,{#2}\,\right\}}
\newcommand{\set}[1]{{\{#1\}}}
\newcommand{\vektor}[1]{{\boldsymbol{#1}}}

\newcommand{\M}{\mathcal{M}}

\newcommand{\UNI}{\mathcal{A}}
\newcommand{\divides}[2]{{#1 \left\lvert {#2} \right.}}
\newcommand{\dividesnot}[2]{{#1 \not \vert #2}}
\newcommand{\ddivides}[2]{{#1 \left\lvert\lvert {#2} \right.}}
\newcommand{\sqbij}{\Phi}

\newcommand{\SQA}{\mathfrak{D}}

\newcommand{\norm}[1]{{ \left\lvert {#1} \right\rvert}}
\newcommand{\tdeg}[1]{{\mathrm{D}(#1)}}
\newcommand{\ord}{\mathrm{N}}
\newcommand{\ann}{\mathrm{ann}}
\newcommand{\supp}{\mathrm{supp}}

\begin{document}

\title[The ring of arithmetical  functions with unitary convolution]{The
  ring of arithmetical 
  functions with unitary convolution: Divisorial and topological properties.}
\author{Jan Snellman}
\address{Department of Mathematics\\
Stockholm University\\
SE-10691 Stockholm,
Sweden}
\email{jans@matematik.su.se}
\date{January 9, 2002}

\begin{abstract}
  We study \((\UNI,+,\oplus)\), the ring of arithmetical functions with unitary
  convolution, 
  giving an isomorphism between \((\UNI,+,\oplus)\) and a generalized power series
  ring on infinitely many variables, similar to the isomorphism
  of Cashwell-Everett\cite{NumThe} 
  between the ring \((\UNI,+,\cdot)\) of arithmetical functions with
  \emph{Dirichlet convolution} 
  and the power series ring \(\C[[x_1,x_2,x_3,\dots]]\) on countably
  many variables. We topologize it with respect to a natural norm, and
  shove that all ideals are quasi-finite. Some 
  elementary results on factorization into atoms are obtained. We
  prove the existence of an abundance of non-associate regular non-units.
\end{abstract}

\maketitle

\begin{section}{Introduction}
  The \emph{ring of arithmetical functions with Dirichlet convolution},
  which we'll denote by \((\UNI,+,\cdot)\),
  is the set of all functions \(\Nat^+ \to \C\), where \(\Nat^+\)
  denotes the positive integers. It is given the structure of a
  commutative \(\C\)-algebra by component-wise addition and
  multiplication by scalars, and by the Dirichlet convolution
  \begin{equation} \label{eq:dirichlet}
    f \cdot g (k) = \sum_{\divides{r}{k}} f(r) g(k/r).
  \end{equation}
  Then, the multiplicative unit is the function \(e_1\) with \(e_1(1)=1\)
  and \(e_1(k)=0\) for \(k>1\), and the additive unit is the zero
  function \(\mathbf{0}\).

  Cashwell-Everett \cite{NumThe} showed that \((\UNI,+,\cdot)\) is a UFD 
  using the isomorphism 
  \begin{equation}\label{eq:ce}
(\UNI,+,\cdot) \simeq \C[[x_1,x_2,x_3,\dots]],
  \end{equation}
  where each \(x_i\) corresponds to the function which is \(1\) on
  the \(i\)'th prime number, and 0 otherwise. 

   Schwab and Silberberg \cite{vring} topologised \((\UNI,+,\cdot)\) by means
   of the norm 
   \begin{equation}
     \label{eq:Snorm}
      \norm{f} = \frac{1}{\min \setsuchas{k}{f(k) \neq 0}}
   \end{equation}
   They noted that this norm is an ultra-metric, and that 
   \(((\UNI,+,\cdot), \norm{\cdot})\) is a valued ring, i.e. that 
   \begin{enumerate}
   \item \(\norm{\mathbf{0}}=0\) and \(\norm{f}>0\) for \(f \neq
     \mathbf{0}\),
   \item \(\norm{f-g} \le \max \set{\norm{f}, \, \norm{g}}\),
   \item \(\norm{fg} = \norm{f} \norm{g}\).
   \end{enumerate}
   They showed that \((\UNI, \norm{\cdot})\) is complete, and that
   each ideal is \emph{quasi-finite}, which means that there exists a
   sequence \((e_k)_{k=1}^\infty\), with
   \(\norm{e_k} \to 0\), such that every
   element in the ideal can be written as a convergent sum
   \(\sum_{k=1}c_k e_k\), with \(c_k \in \UNI\).

  In this article, we treat instead \((\UNI,+,\oplus)\), the ring of
  all arithmetical 
  functions with unitary convolution. This ring has been studied by
  several authors, such as  Vaidyanathaswamy \cite{vaidy},  Cohen
  \cite{UniDiv}, and Yocom \cite{Yocom:totmult}.

  We topologise \(\UNI\) in the same way as Schwab and Silberberg
  \cite{vring},  so that \((\UNI,+,\oplus)\) becomes a
  normed ring (but, in contrast to \((\UNI,+,\cdot)\), not a valued ring). We
  show that all ideals in \((\UNI,+,\oplus)\) are quasi-finite.

  We show that \((\UNI,+,\oplus)\) is
  isomorphic to a monomial quotient of a power series ring on
  countably many variables. It is présimplifiable and atomic, and there is
  a bound on the lengths of factorizations of a given element.
  We give a sufficient condition for  nilpotency, and
  prove the existence of plenty of regular non-units. 

  Finally, we show that the set of arithmetical functions supported on
  square-free integers is a retract of  \((\UNI,+,\oplus)\).

\end{section}

\begin{section}{The ring of arithmetical functions with unitary convolution}
Let \(p_i\) denote
the \(i\)'th prime number, and denote by   \(\mathcal{P}\)
the set of prime numbers. Let \(\mathcal{PP}\) denote the set of
prime powers. Let \(\omega(r)\)
  denote the number of distinct prime factors of \(r\), with
  \(\omega(1)=0\).

  \begin{definition}
    If \(k,m\) are positive integers, we define their \emph{unitary product} as
    \begin{equation}
      \label{eq:unitmm}
      k \oplus m = 
      \begin{cases}
        km & \gcd(k,m)=1 \\
        0 & \text{ otherwise}
      \end{cases}
    \end{equation}
    If \(k \oplus m = p\), then we write \(\ddivides{k}{p}\) and say
    that \(k\) is a \emph{unitary divisor} of \(p\).
  \end{definition}

%  If we let \(\le\) denote the ordinary total order on \(\Nat^+\), then
%  \((\Nat^+, \oplus)\) is a strictly ordered monoid-with-zero, i.e. if \(s <
%  s'\) and \(s \oplus t \neq 0\), \(s' \oplus t \neq 0\),  then \(s
%  \oplus t < s' \oplus t\).  Following Ribenboim
%  \cite{Rib:GP,Rib:Nilpotent} we will study the set \(\C^{(\Nat^+)}\)
%  of all mappings from \(\Nat^+\) to \(\C\), endowed with the
%  \emph{convolution product}. Ribenboim's construction of the generalised
%  power series ring of a strictly ordered monoid consists of choosing
%  those mappings whose support are \emph{artinian} and \emph{narrow},
%  but since every subset of \(\Nat^+\) have those properties (there
%  are no infinite descending chains, or infinite anti-chains), we will
%  work with the set of all mappings. The fact that we are dealing with
%  a monoid-with-zero rather than with a monoid is no major obstacle.

The so-called \emph{unitary convolution} was introduced by
Vaidyanathaswamy \cite{vaidy}, and was further studied Eckford Cohen
\cite{UniDiv}.   
  
  \begin{definition}
    \(\UNI = \set{f: \Nat^+ \to \C}\), the set of complex-valued
    functions on the positive integers. We define the \emph{unitary
      convolution} of \(f,g \in \UNI\) as 
    \begin{equation}
      \label{eq:unitaryconv}
     ( f\oplus g) (n) = \sum_{\substack{m \oplus p = n\\m,n \ge 1}}
     f(m) g(n) = \sum_{\ddivides{d}{n}} f(d) g(n/d)
    \end{equation}
    and the addition as 
\begin{displaymath}
  (f+g)(n) = f(n)+g(n)
\end{displaymath}
The ring \((\UNI, + , \oplus)\) is called \emph{the ring of arithmetic
  functions} with unitary convolution. 
%In this article\footnote{This
%  is at odds with the notation in 
%Sivaramakrishnans treatise \cite{Afunc}, where \(fg\) denotes
%``ordinary multiplication'' \(fg(r)=f(r)g(r)\). We will
%mention were our notations differ from  \cite{Afunc}.
%}, 
%we will often write \(fg\) in place of \(f \oplus g\). 
  \end{definition}

  \begin{definition}
    For each positive integer \(k\), we define \(e_k \in \UNI\) by
    \begin{equation}
      \label{eq:e}
      e_k(n) = 
      \begin{cases}
        1 & k=n \\
        0 & k \neq n
      \end{cases}
    \end{equation}
    We also define\footnote{In \cite{Afunc}, \(\mathbf{1}\) is
      denoted \(e\), and \(e_1\) denoted \(e_0\).} \(\mathbf{0}\) as
    the zero function, and 
    \(\mathbf{1}\) as the function which is constantly 1. 
  \end{definition}

  \begin{lemma}\label{lemma:rels}
    \(\mathbf{0}\) is the additive unit of \(\UNI\), and \(e_1\) is
    the multiplicative unit.  
    We have that
    \begin{equation}
      \label{eq:ee}
      (e_{k_1} \oplus e_{k_2} \oplus \cdots \oplus e_{k_r} ) (n) = 
      \begin{cases}
        1 & n = k_1 k_2 \cdots k_r \text{ and } \gcd(k_i,k_j) = 1
        \text{ for } i \neq j \\  
        0 & \text{ otherwise}
      \end{cases}
    \end{equation}

    hence
    \begin{equation}
      \label{eq:ee2}
      e_{k_1} \oplus e_{k_2} \oplus \cdots \oplus e_{k_r}  = 
      \begin{cases}
        e_{k_1 k_2 \cdots k_r} & \text{ if } \gcd(k_i,k_j) = 1
        \text{ for } i \neq j \\  
        0 & \text{ otherwise}
      \end{cases}
    \end{equation}

  \end{lemma}
  \begin{proof}
    The first assertions are trivial.
    We have \cite{Afunc} that for \(f_1,\dots,f_r \in \UNI\), 
    \begin{equation}
      \label{eq:dist}
      (f_1 \oplus \cdots f_r)(n) = \sum_{a_1 \oplus \cdots a_r=n}
      f_1(a_1) \cdots f_r(a_r)
    \end{equation}
    Since
    \begin{displaymath}
      e_{k_1}(a_1)  e_{k_2}(a_2)  \cdots  e_{k_r})(a_r) = 1 \text{ iff
        } \forall i: k_i=a_i,
    \end{displaymath}
    \eqref{eq:ee} follows. 
  \end{proof}

    \begin{lemma}\label{lemma:uniq}
      Any \(e_n\) can be uniquely expressed as a square-free monomial
      in \(\setsuchas{e_k}{k \in \mathcal{PP}}\).
    \end{lemma}
    \begin{proof}
      By unique factorization, there is a unique way of writing
      \(n=p_{i_1}^{a_1} \cdots p_{i_r}^{a_r}\), and \eqref{eq:ee2} gives 
      that 
      \begin{displaymath}
        e_n = e_{p_{i_1}^{a_1} \cdots p_{i_r}^{a_r}} =
        e_{p_{i_1}^{a_1}} \oplus \cdots e_{p_{i_r}^{a_r}}.
      \end{displaymath}
    \end{proof}

  \begin{theorem}
    \((\UNI, + , \oplus)\) is a quasi-local, non-noetherian
    commutative ring having 
    divisors of zero.  The units
    \(U(\UNI)\) consists of those \(f\) such that \(f(1) \neq 0\).
  \end{theorem}
  \begin{proof}
    It is shown in \cite{Afunc} that  \((\UNI, + , \oplus)\) is a
    commutative ring, having zero-divisors, and that the units
    consists of those \(f\) such that \(f(1) \neq 0\). 
    If \(f(1) = 0\)
    then 
    \begin{displaymath}
      (f \oplus g)(1) = f(1)g(1) = 0.
    \end{displaymath}
    Hence the non-units 
    form an ideal \(\mathfrak{m}\), which is then the unique maximal
    ideal. 
    
    We will show (Lemma~\ref{lemma:mnotfg}) that \(\mathfrak{m}\)
    contains an ideal (the ideal generated by all \(e_k\), for
    \(k>1\)) which is
    not finitely generated, so \(\UNI\) is non-noetherian.
  \end{proof}

\end{section}

\begin{section}{A topology on \(\UNI\)}

The results of this section are inspired by
\cite{vring}, were the authors studied the ring of arithmetical
functions under Dirichlet convolution. We'll use the notations of
\cite{NonArch}. We regard \(\C\) as trivially normed.

  \begin{definition}\label{def:metric}
    Let \(f \in \UNI \setminus \set{\mathbf{0}}\).
    We define the \emph{support} of \(f\) as
    \begin{equation}
      \label{eq:support}
      \supp(f) = \setsuchas{n \in \Nat^+}{f(n) \neq 0}
    \end{equation}
    We define the \emph{order}\footnote{In \cite{Afunc} the term
      \emph{norm} is used.} of a non-zero element  by
    \begin{equation}
      \label{eq:height}
      \ord(f) = \min \supp(f)
    \end{equation}
    We also define the \emph{norm} of \(f\) as
    \begin{equation}
      \label{eq:norm}
      \norm{f} = \ord(f)^{-1}
    \end{equation}
    and the \emph{degree} as
    \begin{equation}
      \label{eq:deg}
      \tdeg{f} = \min \setsuchas{\omega(k)}{k \in \supp(f)}
    \end{equation}
    By definition, the zero element has order infinity, norm 0, and
    degree -1.
  \end{definition}

  \begin{lemma}
    The value semigroup of \((\UNI, \norm{\cdot})\) is 
    \begin{math}
      \norm{\UNI \setminus \set{\mathbf{0}}} = \setsuchas{1/k}{k \in \Nat^+},
    \end{math}
    a discrete subset of \(\R^+\).
  \end{lemma}

%  \begin{proof}
%    \begin{enumerate}
%    \item  Let \(k < i\). Then \((f-g)(k) = f(k) - g(k)
%      = 0\), so \(\ord(f-g) \ge i\).
%    \item  It is clear that \(cf(i) \neq 0\),
%      but that \(cf(k) = 0\), so \(\ord(cf) =  \ord(f)\).
%    \item      \(f\) is a non-unit iff \(f(1) \neq 0\), i.e. iff
%      \(\ord(f)>1\). 
%    \item This is proved in \cite{Afunc}.
%    \item Follows from the previous assertion together with the next.
%    \item     Let \(\ell < ij\), and let \(uv=\ell\), with
%      \(\gcd(u,v)=1\). Then  
%    either \(u < i\) or \(v < j\), hence \(f(u)g(v) =0\). This shows 
%    that \(\ord(fg) \ge   \ord(f)\ord(g)\). Since \(\UNI\) has
%    zero-divisors, strict inequality may occur: 
%    \begin{displaymath}
%\ord(e_2^2) =
%    \ord(\mathbf{0}) = \infty > \ord(e_2)^2 = 2^2.
%    \end{displaymath}
%  \item      If \(f,g\) are non-units, then 
%    \begin{displaymath}
%\ord(fg) \ge \ord(f) \ord(g) >
%     \max \set{\ord(f), \ord(g)}      
%    \end{displaymath}
% since \(\ord(f), \ord(g) \ge
%     2\). If \(f,g\) are both units, then so is \(fg\), hence
%     \begin{displaymath}
%     1=\ord(fg)=\ord(f)=\ord(g).       
%     \end{displaymath}
% If \(f\) is a unit, but \(g\) is
%     not, then \(f(1) \neq 0\), so putting 
%     \(j=\ord(g)\) as before, we have that \(j>1\), and that \((fg)(j)
%     = f(1)g(j) \neq 0\), but that \((fg)(k)=0\) for \(k < i\), hence
%     that \(\ord(fg)=i\).
%    \end{enumerate}
%  \end{proof}

  \begin{lemma}\label{lemma:isval}
    Let \(f,g \in \UNI \setminus \set{\mathbf{0}}\).
    Let \(\ord(f)=i\), \(\ord(g)=j\), so that \(f(i) \neq 0\) but
    \(f(k)=0\) for all \(k < i\), and similarly for \(g\). We assume
    that \(i \le j\).
    Then, the following hold:
    \begin{enumerate}[(i)]
    \item \label{en:ordplus}
      \(\ord(f-g) \ge \min \set{\ord(f), \ord(g)}\).
    \item \label{en:ordconst}
      \(\ord(cf) =  \ord(f)\) for \(c \in \C \setminus \set{0}\).
    \item \label{en:ordunit}
      \(\ord(f)=1\) iff \(f\) is a unit.
    \item \label{en:orddir}
      \(\ord(f \cdot g) = \ord(f) \ord(g) \le \ord(f\oplus g)\), with equality iff \((i,j)=1\).
    \item \label{en:ordun}
      \(\ord(f \oplus g) \ge   \max \set{\ord(f), \ord(g)}\), 
      with strict inequality iff both \(f\) and \(g\)
      are non-units.
  \item \label{en:tdp}
    \(\tdeg{f+g} \ge \min \tdeg{f}, \, \tdeg{g}\).
  \item \label{en:tdm}
    \(\tdeg{f\cdot g} = \tdeg{f}+ \tdeg{g}\).
  \item \label{en:un} \(\tdeg{f}=0\) if and only if \(f\) is a unit.
  \item  \label{en:tdum}
    Suppose that \(f \oplus g  \neq \mathbf{0}\). Then
   \begin{displaymath}
   \tdeg{f\oplus g} \ge \tdeg{f} + \tdeg{g} \ge \max \set{\tdeg{f}, \, 
     \tdeg{g}}.
   \end{displaymath}
   with \(\tdeg{f} + \tdeg{g} > \max \set{\tdeg{f}, \,    \tdeg{g}}\)
   if \(f,g\) are non-units.
  \end{enumerate}
\end{lemma}
\begin{proof}
  \eqref{en:ordplus}, \eqref{en:ordconst}, and \eqref{en:ordunit} are
  trivial, and \eqref{en:orddir} is proved in \cite{Afunc}.
  \eqref{en:tdp}, \eqref{en:tdm}, and \eqref{en:un} are proved in \cite{dvar}.
  Let \(m\) be a monomial in the support of \(f\) such that
  \(\tdeg{m}=\tdeg{f}\), and 
  let \(n\) be a monomial in the support of \(g\) such that
  \(\tdeg{n}=\tdeg{g}\). For any \(a\) in the support of \(f\) and any
  \(q\) in the support of \(g\), such that \(a \oplus q \neq 0\),
  we have that 
  \begin{displaymath}
\tdeg{a \oplus q} = \tdeg{a} + \tdeg{q} \ge \tdeg{f} +
  \tdeg{g}.    
  \end{displaymath}
 This proves \eqref{en:tdum}. \eqref{en:ordun} is proved similarly.
\end{proof}

  \begin{cor}\label{cor:normprod}
    \(\norm{f \oplus g} \le \norm{f} \norm{g} = \norm{f \cdot g}\).
  \end{cor}

  \begin{prop}\label{prop:Cauchy}
   \(\norm{\cdot}\) is an  ultrametric function on \(\UNI\), making
   \((\UNI,+,\oplus)\) a normed ring, as well as a
   faithfully normed, b-separable complete
   vector space over \(\C\).
  \end{prop}
  \begin{proof}
    \(((\UNI,+,\cdot), \norm{\cdot})\) is a
    valuated ring, and a
     faithfully normed  complete vector space over \(\C\)
     \cite{vring}. 
     It is also separable with
    respect to bounded maps \cite[Corollary 2.2.3]{NonArch}. 
    So  \((\UNI,+)\) is a  normed group, hence
    Corollary~\ref{cor:normprod}
    shows that \((\UNI,+,\oplus)\) is a normed ring.
  \end{proof}

Note that, unlike \(((\UNI,+,\cdot),  \norm{\cdot})\), the normed ring
\(((\UNI,+,\oplus), \norm{\cdot})\) is not a valued ring, since 
\begin{displaymath}
  \norm{e_2\oplus e_2} = \norm{\mathbf{0}} = 0 < \norm{e_2}^2 = 1/4.
\end{displaymath}
In fact, we have that 

\begin{lemma}\label{lemma:normsq}
  If \(f\) is a unit, then \(1 = \norm{f^n} = \norm{f}^n\) for all
  positive integers \(n\). If \(n\) is a non-unit, then \(\norm{f^n} <
  \norm{f}^n\) for all \(n>1\).
\end{lemma}
\begin{proof}
  The first assertion is trivial, so suppose that \(f\) is a
  non-unit. From Corollary~\ref{cor:normprod} we 
  have that \(\norm{f^n} \le  \norm{f}^n\). If \(\norm{f}=1/k\),
  \(k>1\), i.e.  \(f(k) \neq 0\) but \(f(j) = 
  0\) for \(j < k\), then \(f^2(k^2) = 0\) since \(\gcd(k,k) = k >
  1\). It follows that \(\norm{f^2} > \norm{f}^2\), from which the
  result follows.
\end{proof}

Recall that in a normed ring, a non-zero element \(f\) is called
\begin{itemize}
\item \emph{topologically nilpotent} if \(f^n \to 0\), 
\item \emph{power-multiplicative} if \(\norm{f^n}=\norm{f}^n\) for all
  \(n\),
\item \emph{multiplicative} if \(\norm{fg}=\norm{f}\norm{g}\) for all
  \(g\) in the ring.
  \end{itemize}

    \begin{theorem}\label{thm:topsim}
      Let \(f \in ((\UNI,+,\oplus),\norm{\cdot})\), \(f \neq\mathbf{0}\). Then
      the following are equivalent: 
      \begin{enumerate}
      \item \(f\) is topologically nilpotent,
      \item \(f\) is not power-multiplicative,
      \item \(f\) is not multiplicative\footnote{This is not the same
        concept as multiplicativity for arithmetical functions,
        i.e. that \(f(nm)=f(n)f(m)\) whenever \((n,m)=1\). However,
        since the latter kind of elements satisfy \(f(1)=1\), they are
        units, and hence multiplicative in the normed-ring sense.} in the normed ring
        \((\UNI,+,\oplus),\norm{\cdot})\), 
      \item \(f\) is a non-unit,
      \item \(\norm{f} < 1\).
      \end{enumerate}
    \end{theorem}
    \begin{proof}
      Using \cite[1.2.2, Prop. 2]{NonArch}, this follows from the
      previous Lemma, and the fact that for a unit 
      \(f\), 
      \begin{displaymath}
        1= \norm{f^{-1}} = \norm{f}^{-1}.
      \end{displaymath}
    \end{proof}

    \begin{subsection}{A Schauder basis for \((\UNI,\norm{\cdot})\)}
    \begin{definition}
      Let \(\UNI'\) denote the subset of \(\UNI\) consisting of
      functions with finite support. We define a pairing
      \begin{equation}
        \label{eq:pairing}
        \begin{split}
          \UNI \times \UNI' & \to \C \\
          \left \langle f,g \right \rangle & = \sum_{k=1}^\infty f(k) g(k)
        \end{split}
      \end{equation}
    \end{definition}

    \begin{theorem}\label{thm:dense}
      The set 
      \begin{math}
        \label{eq:genset}
        \setsuchas{e_{k}}{k \in \Nat^+}
      \end{math}
      is an ordered orthogonal Schauder base in the normed vector space
      \((\UNI, \norm{\cdot})\). In other words,
      if \(f \in \UNI\) then 
      \begin{equation}
        \label{eq:fis}
        f = \sum_{k=1}^\infty c_k e_k, \qquad c_k \in \C
      \end{equation}
      where 
      \begin{enumerate}[(i)]
      \item \(\norm{e_k} \to 0\),
      \item the infinite sum \eqref{eq:fis} converges w.r.t. the
        ultrametric topology,
      \item the coefficients \(c_k\) are
      uniquely determined by the fact that 
      \begin{equation}
        \label{eq:langle}
        \left\langle f, e_k \right\rangle = f(k) = c_k
      \end{equation}
    \item 
      \begin{equation}
        \label{eq:orthogonal}
        \max_{k \in \Nat^+} \set{\norm{c_k}\norm{e_k}} = \norm{
          \sum_{k=1}^\infty c_k e_k} 
      \end{equation}
      \end{enumerate}
      
      The set \(\setsuchas{e_p}{p \in \mathcal{PP}}\) generates a dense
      subalgebra of \(((\UNI,+,\oplus), \norm{\cdot})\).
    \end{theorem}
    \begin{proof}
      It is proved in \cite{vring} that this set
      is a Schauder base in the topological vector space \((\UNI,
      \norm{\cdot})\). It also follows 
      from \cite{vring} that the coefficients \(c_k\) in
      \eqref{eq:genset} are given by \(c_k=f(k)\).

      It remains to prove orthogonality. With the above notation,
      \begin{displaymath}
      \norm{f} = \norm{\sum_{k=1}^\infty c_k e_k}=1/j,        
      \end{displaymath}
 where \(j\)
      is the smallest \(k\) such that \(c_k \neq 0\). Recalling that
      \(\C\) is trivially normed, we have that 
      \begin{displaymath}
        \norm{c_k}\norm{e_k} = 
        \begin{cases}
          \norm{e_k} = 1/k & \text{ if }c_k \neq 0 \\
          0 & \text{ if } c_k = 0
        \end{cases},
      \end{displaymath}
      so \(\max_{k \in \Nat^+} \set{\norm{c_k}\norm{e_k}} = 1/j\),
      with \(j\) as above, so
      \eqref{eq:orthogonal} holds.

      By Lemma~\ref{lemma:uniq} any \(e_k\) can be written as a
      square-free monomial in the elements of 
      \(\setsuchas{e_p}{p \in \mathcal{PP}}\).  The set
      \(\setsuchas{e_k}{k \in \Nat^+}\) is dense in \(\UNI\), so
      \(\setsuchas{e_p}{p \in \mathcal{PP}}\) generates a dense subalgebra.
    \end{proof}

    Let \(J \subset \mathfrak{m}\) denote the ideal generated by all \(e_k\),
      \(k>1\).

    \begin{lemma}\label{lemma:mnotfg}
      \(J\) is not finitely generated.
    \end{lemma}
    \begin{proof}
    If \(J\) is  finitely generated, then
    there is an \(N\) such that 
    \begin{displaymath}
      J = \left( e_2, \dots, e_N  \right). 
    \end{displaymath}    
    Let \(L\) be a prime number, \(L > N\).
    Since \(e_L \in J\), we have that
    \begin{displaymath}
      e_L = \sum_{k=2}^N f_k \oplus e_k, \qquad f_k \in \UNI.
    \end{displaymath}
    We write \(f_k = \sum_{i=1}^\infty c_{ki} e_i\), so that 
    \begin{displaymath}
      e_L = \sum_{k=2}^N e_k \oplus \sum_{i=1}^\infty c_{ik} e_i = 
      \sum_{k=2}^N \sum_{i=1}^\infty c_{ik} e_i \oplus e_k =
      \sum_{k=2}^N \sum_{\gcd(i,k) = 1} c_{ik} e_{ik}.
    \end{displaymath}
    But this is impossible, because we can not write \(L=ik\) with
    \(\gcd(i,k)=1\) and \(2 \le i \le N < L\).
    \end{proof}

    \begin{definition}
      An ideal \(I \subset \UNI\) is called quasi-finite if there
      exists a sequence \((g_k)_{k=1}^\infty\) in \(I\) such that
      \(\norm{g_k} \to 0\) and  such that every element \(f \in I\) can be
      written (not necessarily uniquely) as a convergent sum 
      \begin{equation}
        \label{eq:qf}
        f = \sum_{k=0}^\infty a_k \oplus g_k, \qquad a_k \in \UNI
      \end{equation}
    \end{definition}
    
    \begin{lemma}
      \(\mathfrak{m}\) is quasi-finite.
    \end{lemma}
    \begin{proof}
      By Theorem~\ref{thm:dense} the set \(\setsuchas{e_k}{k > 1}\) is a
      quasi-finite generating set for  \(\mathfrak{m}\).
    \end{proof}
    Since all ideals are contained in \(\mathfrak{m}\), it follows
    that any ideal containing \(\setsuchas{e_k}{k > 1}\) is
    quasi-finite. Furthermore, such an ideal has  \(\mathfrak{m}\) as
    its closure. In particular, \(J\) is quasi-finite, but not closed.

    \begin{theorem}
      All (non-zero) ideals in \(\UNI\) are quasi-finite.
      In fact, given any subspace \(I\) if we can find 
      \begin{equation}
        \label{eq:G}
        G(I) := (g_k)_{k=1}^\infty
      \end{equation}
      such that for all \(f \in I\), 
      \begin{equation}
        \label{eq:qua}
        \exists c_1,c_2,c_3,\dots \in \C, \qquad f = \sum_{i=1}^\infty c_i g_i.
      \end{equation}
      So all subspaces possesses a Schauder basis.
    \end{theorem}
    \begin{proof}
      We construct \(G(I)\) in the following way: for each 
      \begin{displaymath}
k \in
      \setsuchas{\ord(f)}{f \in I  \setminus \set{\mathbf{0}}} =: N(I)        
      \end{displaymath}
      we
      choose a \(g_k \in I\) with   \(\ord(g_k)=k\), and with
      \(g_k(k)=1\). In other words, we make sure that the ``leading
      coefficient'' is 1; this can always be achieved since the
      coefficients lie in a field. For \(k \not \in N(I)\) we put
      \(g_k = \mathbf{0}\).

      To show that this choice of elements satisfy \eqref{eq:qua},
       take any  \(f \in I\), and put
      \(f_0=f\). Then define recursively, as long as \(f_i \neq \mathbf{0}\),
      \begin{align*}
        n_i &:= N(f_i) \\
        \C \ni a_i &:= f_i(n_i)\\
        \UNI \ni f_{i+1} &:= f_i - a_i g_{n_i} 
      \end{align*}
      Of course, if \(f_i= \mathbf{0}\), then we have expressed \(f\)
      as a linear combination of \[g_{n_1},\dots,g_{n_{i-1}},\] and we are
      done. Otherwise, note that by induction \(f_i \in I\), so \(n_i
      \in N(I)\), hence \(g_{n_i} \neq 0\).
      Thus \(\ord(f_{i+1}) > \ord(f_i)\), so \(\norm{f_{i+1}} <
        \norm{f_i}\), whence 
        \begin{displaymath}
          \norm{f_0} > \norm{f_1} > \norm{f_2} > \cdots  \to 0.
        \end{displaymath}
       But \[f_{i+1} = f-\sum_{j=1}^i a_j g_{n_j},\] so 
       \begin{displaymath}
      F_i  :=  \sum_{j=1}^i a_j g_{n_j}  \to  f,   
       \end{displaymath}
       which shows that
      \(\sum_{j=1}^\infty a_j g_j=f\). 
    \end{proof}
  
\end{subsection}

\end{section}

\begin{section}{A fundamental isomorphism}
\begin{subsection}{The monoid of separated monomials}
Let
\begin{equation}
  \label{eq:Y}
  Y = \setsuchas{y_i^{(j)}}{i,j \in \Nat^+}
\end{equation}
be an infinite set of variables, in bijective correspondence with the
integer lattice points in the first quadrant minus the axes. We call
the subset 
\begin{equation}
  \label{eq:column}
  Y_{i} = \setsuchas{y_i^{(j)}}{j \in \Nat^+}
\end{equation}
the \emph{\(i\)'th column} of \(Y\).

Let \([Y]\) denote the free
abelian monoid on \(Y\), and let 
%\([Y']\) denote the subset of
%square-free monomials. We regard  \([Y']\) as a
%  monoid-with-zero, so that the multiplication is given by
%  \begin{equation}
%    \label{eq:zeromult}
%    m \times m' = 
%    \begin{cases}
%      mm' & mm' \in [Y'] \\
%      0 & \text{ otherwise}
%    \end{cases}
%  \end{equation}
%
%
%  Similarly, let 
%
\(\M\) be the subset of \emph{separated monomials},
  i.e. monomials in which no two occurring variables come from the same
  column:
  \begin{equation}
    \label{eq:MM}
    \M = \setsuchas{y_{i_1}^{(j_1)} y_{i_2}^{(j_1)} \cdots y_{i_r}^{(j_r)}}{ 1 \le i_i < i_2
      < \cdots i_r} 
  \end{equation}

  We regard  \(\M\) as a
  monoid-with-zero, so that the multiplication is given by
  \begin{equation}
    \label{eq:mmmult}
    m \oplus m' = 
    \begin{cases}
      mm' & mm' \in \M \\
      0 & \text{ otherwise}
    \end{cases}
  \end{equation}
  Note that the zero is exterior to \(\M\), i.e. \(0 \not \in
  \M\). The set \(\M \cup \set{0}\) is a (non-cancellative) monoid if
  we define \(m \oplus 0 = 0\) for all \(m \in \M\).

Recall that \(\mathcal{PP}\) denotes the set of
prime powers. It follows from the fundamental theorem of arithmetic that any
positive integer
\(n\) can be uniquely written as a \emph{square-free} product of prime
powers. Hence we have that
\begin{equation}
  \label{eq:phi1}
  \begin{split}
  \sqbij: Y & \to \mathcal{PP} \\
  y_i^{(j)} & \mapsto p_i^j
  \end{split}
\end{equation}
is a bijection which can be extended to a bijection 
\begin{equation}
  \label{eq:phi2}
  \begin{split}
  \sqbij: \M & \to \Nat^+ \\
  1 & \mapsto 1\\
  y_{i_1}^{(j_1)} \cdots y_{i_r}^{(j_r)} & \mapsto p_{i_1}^{j_1} \cdots
  p_{i_r}^{j_r} 
\end{split}
\end{equation}
If we regard \(\Nat^+\) as a monoid-with-zero with the operation
\(\oplus\) of \eqref{eq:unitmm}, then \eqref{eq:phi2} is a
monoid-with-zero isomorphism.

%and to a surjection
%\begin{equation}
%  \label{eq:phi3}
%  \begin{split}
%  \sqbij: [Y] & \to \Nat^+ \\
%  1 & \mapsto 1\\
%  y_{i_1}^{(j_1)} \cdots y_{i_r}^{(j_r)} & \mapsto p_{i_1}^{j_1} \cdots
%  p_{i_r}^{j_r} 
%\end{split}
%\end{equation}

\end{subsection}

\begin{subsection}{The ring \(\UNI\) as a generalized power series
    ring, and as a quotient of \(\C[[Y]]\)}

Let \(R\) be the large power series ring on \([Y]\), i.e. \(R=C[[Y]]\)
consists of all formal power series \(\sum c_\vektor{\alpha}
\vektor{y}^\vektor{\alpha}\), where the sum is over all multi-sets
\(\vektor{\alpha}\) on \(Y\). 
%Let \(R_2= \frac{R}{\SQ}\),
%where \(\SQ\) is the ideal generated by \(\setsuchas{y_i^{(j)}^2}{i,j
%  \in \Nat^+}\). 
%Then \(\SQ \subset \SQA\) and hence \(R_3\) is a
%quotient of \(R_2\). 

%  \(R_2\) can be regarded as a generalised monoid-with-zero ring
%  of \([Y']\), i.e. 
%  \(R_2\) is the set of all formal power series 
%  \begin{equation}
%    \label{eq:formal}
%    \sum_{m \in [Y']} f(m) m
%  \end{equation}
%  with multiplication
%  \begin{equation}
%    \label{eq:R2conv2}
%    \begin{split}
%    \left( \sum_{m \in [Y']} f(m) m \right) \times
%    \left( \sum_{m \in [Y']} g(m) m \right) &=
%    \left( \sum_{m \in [Y']} h(m) m \right) \\
%     h(m) &=  f \times g(m) =
%    \sum_{s \times t=m} f(s)g(t) 
%    \end{split}
%  \end{equation}

Let \(S\)  be the  generalized monoid-with-zero ring on \(\M\).
By this, we mean that \(S\) is the set of all formal power series 
\begin{equation}
    \label{eq:formalM}
    \sum_{m \in \M} f(m) m
  \end{equation}
  with component-wise addition, and 
   with multiplication
  \begin{equation}
    \label{eq:R3conv2}
    \begin{split}
    \left( \sum_{m \in \M} f(m) m \right) \oplus
    \left( \sum_{m \in \M} g(m) m \right) &=
    \left( \sum_{m \in \M} h(m) m \right) \\
     h(m) &=  (f \oplus g)(m) =
    \sum_{s \oplus t = m} f(s)g(t)
    \end{split}
  \end{equation}

Define  
\begin{align}
  \supp(\sum_{m \in [Y]} c_m m) &= \setsuchas{m \in Y}{c_m \neq 0} \\
  \supp(\sum_{m \in \M} c_m m) &= \setsuchas{m \in \M}{c_m \neq 0} \\
\end{align}

Let furthermore  
\begin{equation}
  \label{eq:sqa}
\SQA = \setsuchas{f \in R}{\supp(f) \cap \M = \emptyset}
\end{equation}

%\begin{definition}\label{def:ealpha}
%  For each multi-index \(\vektor{\alpha}=(\alpha_1,\dots,\alpha_n)\),
%  we define \(e_{\vektor{\alpha}} = e_{p_1^{\alpha_1} \cdots
%    p_n^{\alpha_n}} \in \UNI\).
%\end{definition}

\begin{theorem}\label{thm:iso}
  \(S\) and \(\frac{R}{\SQA}\) and \(\UNI\) are isomorphic as
  \(\C\)-algebras.    
\end{theorem}
\begin{proof}
  The bijection \eqref{eq:phi2} induces a bijection between \(S\)
  and \(\UNI\) which is an isomorphism because of the way
  multiplication is defined on \(S\). In detail, the isomorphism is
  defined by 
  \begin{equation}
    \begin{split}
    \label{eq:triviso}
  S \ni  \sum_{m \in \M} c_m m  & \mapsto f \in \UNI \\
  f(\Phi(m)) &= c_m
    \end{split}
  \end{equation}

  For the second part,  consider the epimorphism
  \begin{align*}
    \phi: R &\to S \\
    \phi \left( \sum_{m \in [Y]} c_m m \right) &=  
    \sum_{m \in \M} c_m m 
  \end{align*}
  Clearly, \(\ker(\phi) = \SQA\), hence \(S \simeq
  \frac{R}{\ker(\phi)} = \frac{R}{\SQA}\).
\end{proof}

Let us exemplify this isomorphism by noting that \(e_n\), where \(n\)
has the square-free factorization
\(n=p_1^{a_1} \cdots  p_r^{a_r}\), corresponds
to the square-free monomial \(y_{1}^{(a_1)} \cdots y_{r}^{(a_r)}\), and that

\begin{equation}
  \label{eq:oneis}
      \mathbf{1} = \sum_{m \in \M} m 
    = \prod_{i=1}^\infty \left( 1+\sum_{j=1}^\infty y_i^{(j)} \right)
\end{equation}

What does its inverse \(\mu^*\) correspond to?

\begin{definition}
  For \(m \in \M\), we denote by \(\tdeg{m}\) the number of occurring
  variables in \(m\) (by definition, \(\tdeg{1}=0\) and
  \(\tdeg{0}=-\infty\)). For 
  \begin{displaymath}
S \ni f 
  = \sum_{m \in \M} c_m m    
  \end{displaymath}
 we  put 
  \begin{equation}
    \label{eq:abs}
    \tdeg{f} = \min \setsuchas{\tdeg{m}}{c_m \neq 0}
  \end{equation}
  Using the isomorphism between \(S\) and \(\UNI\), we define
  \(\tdeg{g}\) for any \(g \in \UNI\) by 
  \begin{displaymath}
    \tdeg{g} = \min \setsuchas{\omega(n)}{f(n) \neq 0}.
  \end{displaymath}
\end{definition}

It is known (see \cite{Afunc}) that 
%\begin{equation}
%  \label{eq:dstar}
%  (\mathbf{1} \oplus \mathbf{1})(r) = d^*(r) = 2^{\omega(r)}
%\end{equation}
%and that 
\begin{equation}
  \label{eq:omegastart}
  \mu^*(r)=(-1)^{\omega(r)}
\end{equation}
We then have that \(\mu^*\) corresponds to
\begin{equation}
  \label{eq:mus}
  {\mathbf{1}}^{-1} = \frac{1}{\prod_{i=1}^\infty \left(
        1+\sum_{j=1}^\infty y_i^{(j)} \right)} = \prod_{i=1}^\infty
    \frac{1}{ 1+\sum_{j=1}^\infty y_i^{(j)}} 
    =\sum_{m \in \M} (-1)^{\tdeg{m}} m 
\end{equation}

Recall that \(f \in \UNI\) is a \emph{multiplicative} arithmetic
function if \(f(nm)=f(n)f(m)\) whenever \((n,m)=1\). Regarding \(f\)
as an element of \(S\) we have that \(f\) is multiplicative if and
only if it can be written as 
\begin{equation}
  \label{eq:multiplicative}
  f = \prod_{i=1}^\infty \left(1 + \sum_{j=1}^\infty c_{i,j} y_i^{(j)}
    \right)
\end{equation}
  It is now easy to see that the multiplicative functions form a group
  under multiplication.
\end{subsection}

\begin{subsection}{The continuous endomorphisms}
In \cite{vring}, Schwab and Silberberg characterized all continuous
endomorphisms of \(\Gamma\). We give the corresponding result for \(\UNI\):
\begin{theorem}
  Every continuous endomorphism \(\theta\) of the \(\C\)-algebra
  \(S \simeq \UNI\) is defined by 
  \begin{equation}
    \label{eq:endo}
    \theta(y_i^{(j)}) = \gamma_{i,j} 
  \end{equation}
  where 
  \begin{equation}
    \label{eq:kernel}
    \gamma_{i,j} \gamma_{i,k}  = 0 \qquad \text{ for all } i,j,k
  \end{equation}
  and
  \begin{equation}
    \label{eq:cont}
    \gamma_{a_1(n),b_1(n)} \cdots \gamma_{a_r(n),b_r(n)}  \to
    0 \qquad \text{ as } n = p_{a_1(n)}^{b_1(n)} \cdots
    p_{a_r(n)}^{b_r(n)} 
    \to \infty
  \end{equation}
\end{theorem}
\begin{proof}
  Recall that \(S \simeq \frac{R}{\SQA}\), where \(R= \C[[Y]]\)
  and \(\SQA\) is the closure of the ideal generated by all non-separated quadratic
  monomials \(y_i^{(j)}y_{i}^{(k)}\). Since the set of square-free monomials
  in the \(y_i^{(j)}\)'s form a Schauder base, any  continuous
  \(C\)-algebra endomorphism \(\theta\) of \(S\) is 
  determined by its values on the \(y_i^{(j)}\)'s, and must fulfill
  \eqref{eq:cont}. Since \(y_i^{(j)}y_{i}^{(k)} = 0\) in \(S\), we must
  have that 
  \begin{displaymath}
    \theta(0) = \theta(y_i^{(j)}y_{i}^{(k)}) =
    \theta(y_i^{(j)})\theta(y_{i}^{(k)}) = 
    \gamma_{i,j} \gamma_{i,k} = 0.
  \end{displaymath}
\end{proof}
  
\end{subsection}

\end{section}

\begin{section}{Nilpotent elements and zero divisors}
\newcommand{\psupp}{\mathrm{psupp}}
\newcommand{\lp}{\mathrm{lp}}
\newcommand{\existsunique}{\exists !}
\begin{definition}
For \(m \in \Nat^+\), define the \emph{prime support} of \(m\) as
\begin{equation}
  \label{eq:primesupp}
  \psupp(m)=\setsuchas{p \in \mathcal{P}}{\divides{p}{m}}
\end{equation}
and (when \(m>1\)) the \emph{leading prime} as 
\begin{equation}
  \label{eq:lp}
  \lp(m) = \min \psupp(m)
\end{equation}
For \(n \in \Nat^+\), put
\begin{equation}
  \label{eq:supportfilt}
  A^{n} = \setsuchas{k \in \Nat^+}{\divides{p_n}{k} \text{ but }
    \dividesnot{p_i}{k} \text{ for } i < n} =
  \setsuchas{k \in \Nat^+}{\lp(k) = p_n}
\end{equation}
Then \(\Nat^+ \setminus \set{1}\) is a disjoint union
\begin{equation}
  \label{eq:disj}
  \Nat^+ \setminus \set{1} = \bigsqcup_{i=1}^\infty A^{i}
\end{equation}
\end{definition}

\begin{definition}\label{def:canonical}
Let \(f \in \UNI\) be a non-unit. The \emph{canonical decomposition}
of \(f\) is the unique way of expressing \(f\) as a convergent sum 
\begin{equation}
  \label{eq:cansum}
  f = \sum_{i=1}^\infty f_i, \quad   f_i = \sum_{k \in A^i}{f(k) e_k}
\end{equation}
 The element \(f\) is said to be of \emph{polynomial type} if
all but finitely many of the \(f_i\)'s are zero. In that case, the
largest \(N\) such that \(f_N \neq \mathbf{0}\) is called the
\emph{filtration degree} of \(f\).
\end{definition}

\begin{lemma}\label{lemma:jonce}
\begin{equation}
  \label{eq:fi}
f_i = \sum_{j=1}^\infty e_{p_i^j}
  \oplus g_{i,j}, \quad r \le i, \quad  \divides{p_r}{n} \,\,
  \implies \, g_{i,j}(n)=0.
\end{equation}
For any \(n\) there is at most one pair \((i,j)\) such that 
\begin{displaymath}
  \left(e_{p_i^j}  \oplus g_{i,j}\right)(n) \neq 0.
\end{displaymath}
  More precisely, if 
  \begin{displaymath}
    n=p_{i_1}^{j_1} \cdots p_{i_r}^{j_r}, \qquad i_1 < \cdots < i_r,
  \end{displaymath}
  then \(\left(e_{p_{i_1}^{j_1}}  \oplus g_{i_1,j_1}\right)(n)\) may
  be non zero.
\end{lemma}

\begin{definition}
For \(k \in \Nat\), define 
\begin{equation}
  \label{eq:Iid}
  I_{k} = \setsuchas{f \in \UNI}{f(n)=0 \text{ for every } n \text{
      such that } (n,p_1p_2 \cdots p_k)=1}
\end{equation}
\end{definition}
\begin{lemma}
  \(I_k\) is an  ideal in \((\UNI,+,\oplus)\).
\end{lemma}
\begin{proof}
It is shown in \cite{dvar} that the \(I_k\)'s form an ascending chain
of ideals in \((\UNI,+,\cdot)\). They are also easily seen to be
ideals in \((\UNI,+,\oplus)\): if 
\[f \in I_k, \,g \in \UNI \text{ and }
(n,p_1p_2 \cdots p_k)=1\] then 
\begin{displaymath}
  (f \oplus g)(n) = \sum_{\ddivides{d}{n}} f(d) g(n/d) = 0,
\end{displaymath}
since \((d,p_1p_2 \cdots p_k)=1\) for any unitary divisor of \(n\).
\end{proof}

\begin{thm}\label{thm:fszN}
  Let \(N \in \Nat^+\), and let \(f \in (\UNI,+,\oplus)\) be a
  non-unit. Then 
  \begin{align*}
    \label{eq:eq:eqideal}
    I_N &= \ann(e_{p_1 \cdots p_N}) \\ 
    &= \set{\mathbf{0}} \cup 
    \setsuchas{f \in \UNI}{f \text{ is of
        polynomial type and has filtration degree } N} \\
    &=\overline{\UNI \setsuchas{e_{p_i^a}}{a ,i \in \Nat^+,\, i \le N}} 
  \end{align*}
  where \(\overline{\UNI W}\) denotes the topological closure of the
  ideal generated by the set \(W\).
\end{thm}
\begin{proof}
  If \(f \in I_N\) then for all \(k\)
  \begin{equation}
     \quad (f \oplus e_{p_1 \cdots p_N}) (k) = 
\sum_{a \oplus p_1  \cdots  p_N = k}  f(a)  e_{p_1  \cdots  p_N} (p_1
\cdots  p_N) =
\sum_{a \oplus p_1  \cdots  p_N = k}  f(a) =0
  \end{equation}
  so \(f \in \ann(e_{p_1 \cdots p_N})\). 
  Conversely, if \(f \in
    \ann(e_{p_1 \cdots p_N})\) then \((f \oplus e_{p_1 \cdots p_N}) (k)
    =0\) for all \(k\),  hence if \((n,p_1 \cdots p_N)=1\) then 
    \begin{equation}
      \label{eq:iszerois}
      0 = (f \oplus e_{p_1 \cdots p_N}) (n p_1 \cdots p_N)
      = f(n) e_{p_1 \cdots p_N}( p_1 \cdots p_N) = f(n)
    \end{equation}
    hence \(f \in I_N\).

    If \(f \in I_N\) then for \(j > N\) we get that \(f_j =
    \mathbf{0}\), since 
    \begin{displaymath}
      f_j(k) = 
      \begin{cases}
        0 & \text{ if } k \not \in A^j\\
        f(k) =0 & \text{ if } k \in A^j
      \end{cases}
    \end{displaymath}
    Hence \(f = \sum_{i=1}^N f_i\). Conversely, if \(f\) can be
    expressed thusly, then \(f(k)=f_{j_1}(k)=0\) for \(k=p_{j_1}^{a_1}
    \cdots p_{j_r}^{a_r}\) with \(N < j_1 < \cdots < j_r\).

    The last equality follows from Theorem~\ref{thm:dense}.
\end{proof}

\begin{thm}\label{thm:fsz}
  Let \(f \in \UNI\) be a non-unit. The following are equivalent:
  \begin{enumerate}[(i)]
  \item \label{en:pt} \(f\) is of polynomial type.
  \item \label{en:filt} \(f \in \cup_{k=0}^\infty I_k\),
  \item \label{en:fin}
    There is a finite subset 
    \(Q  \subset \mathcal{P}\) such that \(f(k) = 0\) for all \(k\)
    relatively prime to all \(p \in Q\). 
  \item \label{en:ann}
    \(f \in \cup_{N=1}^\infty \ann(e_{p_1p_2\cdots p_N})\). 
  \item \label{en:clos}
    \(f\) is contained in the topological closure of the ideal
    generated by the set 
    \(\setsuchas{e_{p_i^a}}{a ,i \in \Nat^+,\, i \le N}\).
  \end{enumerate}
  If \(f\) has finite support, then it is of polynomial type.
  If \(f\) is of polynomial type, then it is nilpotent. 
\end{thm}
\begin{proof}
  Clearly, a finitely supported \(f\) is of polynomial type.
  The equivalence \eqref{en:pt} \(\iff\) \eqref{en:filt}
  \(\iff\)\eqref{en:fin} \(\iff\) \eqref{en:ann} \(\iff\) \eqref{en:clos}
  follows from the previous theorem.

  If \(f\) is of polynomial type, say of filtration degree \(N\),
  then 
  \begin{equation}
    \label{eq:sumform}
    f = \sum_{i=1}^N f_i
  \end{equation}
  and we see that if \(f^{N+1}\) is the \(N+1\)'st unitary
  power of \(f\), then \(f^{N+1}\) is the linear combination of monomials
  in the \(f_i\)'s, and none of these monomials are
  square-free. Since \(f_i \oplus f_i = \mathbf{0}\) for all \(i\),
  we have that \(f^{N+1}=\mathbf{0}\). So \(f\) is nilpotent.
\end{proof}

\begin{lemma}\label{lemma:poltypid}
  The elements of polynomial type forms an ideal.
\end{lemma}
\begin{proof}
  By the previous theorem, this set can be expressed as 
  \begin{displaymath}
    \bigcup_{n=1}^\infty I_n,
  \end{displaymath}
  which is an ideal since each \(I_n\) is.
\end{proof}

\begin{question}
  Are all [nilpotent elements, zero divisors] of  polynomial type? 
If one could prove that the zero divisors are precisely the elements
of polynomial type, then by Lemma~\ref{lemma:poltypid} it would follow
that 
\begin{math}
  Z(\UNI) 
\end{math}
is an ideal, and moreover a prime ideal, since the product of two
regular elements is regular (in any commutative ring). Then one could
conclude \cite{zerodiv} that 
\((\UNI,+,\oplus)\) has \emph{few zero divisors}, hence is
\emph{additively regular}, hence is a \emph{Marot ring}. 
\end{question}

\begin{theorem}\label{thm:regN}
  \((\UNI,+,\oplus)\) contains infinitely many non-associate regular
  non-units. 
\end{theorem}
\begin{proof}
  {\bf Step 1}. We first show that there is at least one such element.
  Let \(f \in \UNI\) denote the arithmetical function
  \begin{displaymath}
    f(k) = 
    \begin{cases}
      1 & k \in \mathcal{PP} \\
      0 & \text{otherwise }
    \end{cases}
  \end{displaymath}
  Then \(f\) is a non-unit, and using a result by Yocom
  \cite{Yocom:totmult, dvar} we have that \(f\) is contained in a
  subring of \((\UNI,+,\oplus)\) which is a  discrete valuation ring
  isomorphic to \(\C[[t]]\), the power series ring in one
  indeterminate. This ring is a domain, so \(f\) is not nilpotent.

  We claim that \(f\) is in fact regular. 
  To show this, suppose that \(g \in  \UNI\), \(f\oplus g =
  \mathbf{0}\). We will show that 
  \(g =  \mathbf{0}\).

  Any positive integer \(m\) can be written \(m=q_1^{a_1} \cdots q_r^{a_r}\),
  where the \(q_i\) are distinct prime numbers.
  If \(r=0\), then \(m=1\), and \(g(1)=0\), since
  \begin{displaymath}
    0=(f\oplus g)(2) = f(2) g(1) = g(1).
  \end{displaymath}
  For the case \(r=1\),  we want to show that \(g(q^a)=0\) for all
  prime numbers \(q\).
  Choose three different prime powers \(q_1^{a_1}\), \(q_2^{a_2}\), and \(q_3^{a_3}\).
  Then
  \begin{displaymath}
    0 =  f\oplus g(q_i^{a_i} q_j^{a_j}) = f(q_i^{a_i}) g(q_j^{a_j}) + f(q_j^{a_j}) g(q_i^{a_i}) =
    g(q_j^{a_j}) + g(q_i^{a_i}),
  \end{displaymath}
  when \(i \neq j\), \(i,j \in \set{1,2,3}\).
  In matrix notation, these three equations can be written as
  \begin{displaymath}
    \begin{bmatrix}
      1    & 1 & 0       \\
      1 & 0       & 1 \\
      0      & 1  & 1
    \end{bmatrix}
    \begin{bmatrix}
      g(q_1^{a_1}) \\ g(q_2^{a_2}) \\ g(q_3^{a_3})
    \end{bmatrix}
    =
    \begin{bmatrix}
      0 \\ 0 \\ 0
    \end{bmatrix}
  \end{displaymath}
  from which we conclude (since the determinant of the coefficient
  matrix is non-zero) that \(0 = g(q_1^{a_1}) = g(q_2^{a_2}) = g(q_3^{a_3})\).

  Now for the general case, \(r>1\). We need to show that
  that 
  \begin{equation}
    \label{eq:mc}
    g(q_1^{a_1} \cdots q_r^{a_r}) = 0
  \end{equation}
  whenever \(q_1^{a_1}, \dots, q_r^{a_r}\) are pair-wise relatively
  prime prime powers.

  Choose \(N\) pair-wise relatively prime prime powers \(q_1^{a_1},
  \dots, q_N^{a_N}\).  
  For each \(r+1\)-subset \(q_{s_1},\dots, q_{s_{r+1}}\) of this set
  we get a
  homogeneous linear equation 
  \begin{multline}
    \label{eq:mc2}
    0 = f\oplus g(q_{s_1} \dots q_{s_{r+1}}) = \\
    g(q_{s_2}\cdots q_{s_{r+1}} ) + g(q_{s_1} q_{s_3} \cdots
    q_{s_{r+1}}  ) + \dots + g(q_{s_1} \cdots g_{s_r} )
  \end{multline}
  The  matrix of the homogeneous linear equation system formed by all
  these equations is the incidence matrix of \(r\)-subsets (of a
  set of \(N\) elements) into \(r+1\)-subsets. It
  has full rank \cite{Wilson:design}. Since it consists of
  \(\binom{N}{r+1}\) equations and \(\binom{N}{r}\) variables,
  we get that 
  for sufficiently large \(N\), the null-space is zero-dimensional,
  thus the homogeneous system has only the trivial solution.
  It follows, in particular,  that \eqref{eq:mc} holds.

  Thus, \(g(m)=0\) for all \(m\), so \(f\) is a regular element.

  {\bf Step 2}. We construct infinitely many different regular non-units.
  Consider  the element \(\tilde{f}\), with 
  \begin{displaymath}
    \tilde{f}(k) = 
    \begin{cases}
      c_k & k \in \mathcal{PP} \\
      0 & \text{otherwise }
    \end{cases}
  \end{displaymath}
  and where the \(c_k\)'s are ``sufficiently generic'' non-zero
  complex numbers, then we claim that \(\tilde{f}\), too, is a
  regular non-unit. 
  With \(g\), \(m\), \(r\) as before, we
  have that, for \(r=0\), 
  \begin{displaymath}
    0 =   f\oplus g(p^a) = = f(p^a) g(1) = c_{p^a} g(1). 
  \end{displaymath}
  We demand that \(c_{p^a} \neq 0\), then \(g(1)=0\).

  For a general \(r\), we argue as follows: the incidence matrices
  that occurred before will be replaced with ``generic'' matrices
  whose elements are \(c_k\)'s or zeroes, and which specialize, when
  setting all \(c_k=1\), to full-rank matrices. They must therefore
  have full rank, and the proof goes through.

  {\bf Step 3}.   Let \(g\) be a  unit in \(\UNI\), and \(\tilde{f}\) as above.
  We claim that
  if \(g \oplus f\) is of the above 
  form, i.e. supported on \(\mathcal{PP}\), then \(g\) must be a
  constant. Hence there are infinitely 
  many non-associate regular non-units of  the above form.

  To prove the claim, we argue exactly as before, using the fact that
  \(g \oplus \tilde{f}\) is supported on \(\mathcal{PP}\).  
  For \(m=q_1^{a_1} \cdots q_r^{a_r}\) as before, the case \(r=0\) yields
  nothing: 
 \begin{displaymath}
    0 = g\oplus \tilde{f}(1) = \tilde{f}(1) g(1) = 0 g(1)=0,
  \end{displaymath}
  neither does the case \(r=1\):
  \begin{displaymath}
    w = g\oplus \tilde{f}(q^a) = \tilde{f}(q^a) g(1),
  \end{displaymath}
  so \(g(1)\) may be non-zero.
  But for \(r=2\) we get
  \begin{displaymath}
    0 = g\oplus \tilde{f}(q_1^{a_1}q_2^{a_2}) = 
    \tilde{f}(q_1^{a_1}) g(q_2^{a_2}) +
    g(q_1^{a_1}) \tilde{f}(q_2^{a_2}),
  \end{displaymath}
  and also 
  \begin{align*}
    0 &= g\oplus \tilde{f}(q_1^{a_1}q_3^{a_3}) = 
    \tilde{f}(q_1^{a_1}) g(q_3^{a_3}) +
    g(q_1^{a_1}) \tilde{f}(q_3^{a_3}) \\
    0 &= g\oplus \tilde{f}(q_2^{a_2}q_3^{a_3}) = 
    \tilde{f}(q_2^{a_2}) g(q_3^{a_3}) +
    g(q_1^{a_1}) \tilde{f}(q_3^{a_3}) 
  \end{align*}
  which means that 
  \begin{displaymath}
    \begin{bmatrix}
      \tilde{f}(q_2^{a_2})    & \tilde{f}(q_1^{a_1}) & 0       \\
      \tilde{f}(q_3^{a_3})         & 0       & \tilde{f}(q_1^{a_1}) \\
      0      & \tilde{f}(q_3^{a_3})  & \tilde{f}(q_2^{a_2})
    \end{bmatrix}
    \begin{bmatrix}
      g(q_1^{a_1}) \\ g(q_2^{a_2}) \\ g(q_3^{a_3})
    \end{bmatrix}
    =
    \begin{bmatrix}
      0 \\ 0 \\ 0
    \end{bmatrix}
  \end{displaymath}
  By our assumptions, the coefficient matrix is non-singular, so
  only the zero solution exists, hence \(g(q_1^{a_1})=0\).

  An analysis similar to
  what we did before shows that \(g(q_1^{a_1}\cdots q_r^{a_r})=0\) for
  \(r > 1\).
\end{proof}

With the same method, one can easily show that the characteristic
function on \(\mathcal{P}\) is regular.


%%% Local Variables: ***
%%% TeX-master:"unitary" ***
%%% End: ***

\end{section}

\begin{section}{Some simple results on factorisation}
  Cashwell-Everett \cite{NumThe} showed that \((\UNI,+,\cdot)\) is a
  UFD. We will briefly treat the factorisation properties of \((\UNI,+,\oplus)\). 
Definitions and facts regarding factorisation in commutative rings with
zero-divisors from the articles by Anderson and Valdes-Leon
\cite{ZeroDiv1,ZeroDiv2} will be used.

First, we note that since \((\UNI,+,\oplus)\) is quasi-local, it is
présimplifiable, i.e.  \(a \neq \mathbf{0}\), \(a = r \oplus a\)
implies that  \(r\) is a unit. It follows that for  \(a,b \in \UNI\),
the following three 
conditions are equivalent: 
\begin{enumerate}
\item \(a,b\) are \emph{associates}, i.e. \(\UNI \oplus a = \UNI
  \oplus b\).
\item \(a,b\) are \emph{strong associates}, i.e. \(a = u \oplus b\)
  for some unit  \(u\). 
\item \(a,b\) are \emph{very strong associates}, i.e. \(\UNI \oplus a
  = \UNI \oplus b\) and
  either \(a=b=\mathbf{0}\), or \(a \neq \mathbf{0}\) and 
  \(a=r\oplus b \implies r \in U(\UNI)\).
\end{enumerate}
We say that \(a \in \UNI\) is \emph{irreducible}, or an \emph{atom}, if
\(a = b\oplus c\) implies 
that \(a\) is associate with either \(b\) or \(c\). 

\begin{theorem}\label{thm:atomic}
  \((\UNI,+,\oplus)\) is \emph{atomic}, i.e. all non-units can be
  written as a product of finitely many atoms. In fact, \((\UNI,+,\oplus)\) is a
\emph{bounded factorial ring} (BFR), i.e. there is a bound on the
length of all factorisations of an element.
\end{theorem}

\begin{proof}
  It follows from Lemma~\ref{lemma:isval} that  the 
  non-unit \(f\) has a factorisation into at most \(\tdeg{f}\) atoms.
\end{proof}

\begin{example}\label{ex:1}
  We have that \(e_2 \oplus (e_{2^k} + e_3) = e_6\) for all
      \(k\), hence \(e_6\) has an infinite number of non-associate
      irreducible divisors, and infinitely many factorisations into atoms.
\end{example}

\begin{example}\label{ex:2}
   The element \(h= e_{30}\) can be factored as \(e_2 \oplus
      e_3 \oplus e_5\), or as \((e_6 + e_{20}) \oplus (e_2 + e_5)\).
\end{example}

These examples show that \((\UNI,+,\oplus)\) is neither a
\emph{half-factorial ring}, nor a  \emph{finite factorisation ring},
nor a \emph{weak finite factorisation ring}, nor an \emph{atomic
  idf-ring}.

\end{section}

\begin{section}{The subring of arithmetical functions supported on
    square-free integers}
\newcommand{\SQF}{\mathcal{SQF}}
\newcommand{\SQFS}{\mathfrak{C}}
\newcommand{\sqfrest}{p}
Let \(\SQF \subset \Nat^+\) denote the set of square-free integers,
and put
\begin{equation}
  \label{eq:sqfree}
\SQFS = \setsuchas{f \in \UNI}{\supp(f) \subset \SQF}
\end{equation}
For any \(f \in \UNI\), denote by \(\sqfrest(f) \in \SQFS\) the restriction of
\(f\) to \(\SQF\).

\begin{theorem}
  \((\SQFS,+,\oplus)\) is a subring of \((\UNI,+,\oplus)\), and a
  closed \(\C\)-subalgebra with respect to the norm
  \(\norm{\cdot}\). The map 
  \begin{equation}
    \label{eq:restmap}
    \begin{split}
    \sqfrest: \UNI & \to \SQFS\\
    f & \mapsto \sqfrest(f)
    \end{split}
  \end{equation}
  is a continuous \(\C\)-algebra epimorphism, and a retraction of the
  inclusion map \(\SQFS \subset \UNI\).
\end{theorem}
\begin{proof}
  Let \(f,g \in \SQFS\). If \(n \in \Nat^+ \setminus \SQF\) then
  \((f+g)(n) = f(n) + g(n) = 0\), and \(cf(n)=0\) for all \(c \in
  \C\). Since \(n \in \Nat^+ \setminus \SQF\), there is at least on
  prime \(p\) such that \(\divides{p^2}{n}\). If \(m\) is a
  unitary divisor of \(m\), then either \(m\) or \(n/m\) is divisible
  by \(p^2\). Thus \[(f \oplus g)(n) = \sum_{\ddivides{m}{n}} f(m)
  g(n/m) =0.\] 

  If \(f_k \to f\) in \(\UNI\), and all \(f_k \in \SQFS\), let \(n \in
  \supp(f)\). Then there is an \(N\) such that \(f(n)=f_k(n)\) for all
  \(k \ge N\). But \(\supp(f_k) \subset \SQF\), so \(n \in
  \SQF\). This shows that \(\SQFS\) is a closed subalgebra of
  \(\UNI\).

  It is clear that \(\sqfrest(f+g)= \sqfrest(f) + \sqfrest(g)\) and
  that \(\sqfrest(cf)= c\sqfrest(f)\) for any \(c \in \C\). If \(n\)
  is not square-free, we have already showed that 
  \begin{displaymath}
    0=(\sqfrest(f) \oplus \sqfrest(g))(n) = \sqfrest((f \oplus g))(n).
  \end{displaymath}
  Suppose therefore that \(n\) is square-free. Then so is all its
  unitary divisors, hence 
  \begin{multline*}
    \sqfrest(f \oplus g)(n) = (f \oplus g)(n) =
    \sum_{\ddivides{m}{n}} f(m) g(n/m) =  \\
    \sum_{\ddivides{m}{n}}
    \sqfrest(f)(m) \sqfrest(g)(n/m) = (\sqfrest(f) \oplus \sqfrest(g))(n).
  \end{multline*}

  We have that \(\sqfrest(f)=f\) if and only if \(f \in \SQFS\), hence
  \(\sqfrest(\sqfrest(f))=\sqfrest(f)\), so \(\sqfrest\) is a
  retraction to the inclusion \(i: \SQFS \to \UNI\). In other words,
  \(\sqfrest \circ i = \mathrm{id}_\SQFS\).
\end{proof}

\begin{corr}
  The multiplicative inverse of an element in \(\SQFS\) lies in \(\SQFS\).
\end{corr}
\begin{proof}
  If \(f \in \SQFS\), \(f\oplus g = e_1\) then 
  \begin{displaymath}
    e_1 = \sqfrest(e_1) = \sqfrest(f\oplus g) = \sqfrest(f)\oplus
    \sqfrest(g) =  f \oplus \sqfrest(g),
  \end{displaymath}
  hence \(g= \sqfrest(g)\), so \(g \in \SQFS\).

  Alternatively, we can reason as follows. If \(f\) is a unit in
  \(\SQFS\) then we can without loss of generality assume that \(f(1)=1\). By
  Theorem~\ref{thm:topsim}, \(g=-f+e_1\) is topologically nilpotent, hence
  by Proposition 1.2.4 of \cite{NonArch} we have that the inverse
  of \(e_1 - g = f\) can be expressed as \(\sum_{i=0}^\infty g^i\). It
  is clear that \(g\), and every power of it, is supported on
  \(\SQF\), hence so is \(f^{-1}\).
\end{proof}

\begin{corr}
  \((\SQFS,+,\oplus)\) is semi-local.
\end{corr}
\begin{proof}
  The units consists of all  \(f \in \SQFS\) with  \(f(1) \neq
  0\), and the non-units form the unique maximal ideal.
\end{proof}

\begin{note}
  More generally, given any subset \(Q \subset \Nat^+\), 
  we get a retract of \((\UNI,+,\oplus)\) when considering those
  arithmetical functions that are supported on the integers
  \(n=p_1^{a_1} \cdots p_r^{a_r}\) with \(a_i \in Q \cup \set{0}\). 
  This property is unique for the unitary convolution, among all
  regular convolutions in the sense of Narkiewicz \cite{Nark:conv}.
  
  In  particular, the set of arithmetical functions supported on the
  exponentially odd integers (those \(n\) for which all \(a_i\) are
  odd) forms a retract of \((\UNI,+,\oplus)\). It follows that the
  inverse of such a function is of the same form.
\end{note}

Let \(T=\C[[x_1,x_2,x_3,\dots]]\), the large power series ring on
countably many variables, and let \(J\) denote the ideal of
elements supported on non square-free monomials.

\begin{theorem}
  \((\SQFS,+,\oplus) \simeq T/J\).
  This algebra can also be described as the generalized power series
  ring on the monoid-with-zero whose elements are all finite subsets of a fixed
  countable set \(X\), with multiplication 
  \begin{equation}
    \label{eq:multimult}
    A \times B = 
    \begin{cases}
      A \cup B & \text{ if } A \cap B = \emptyset\\
      0 & \text{ otherwise }.
    \end{cases}
  \end{equation}
\end{theorem}
\begin{proof}
  Define \(\eta\) by
  \begin{equation}
    \label{eq:sqdef}
    \begin{split}
    \eta: T & \to \UNI \\
    \eta(\sum_{m} c_m m) &= \sum_{m \text{ square-free}} c_m e_m,
    \end{split}
  \end{equation}
  where for a square-free monomial \(m=m_{i_1} \cdots m_{i_r}\) with
  \(1 \le i_1 < \cdots < i_r\) we put \(e_m = e_{p_{i_1} \cdots p_{i_r}}\).
  Then  \(\eta(T) = \SQFS\), \(\ker \eta = J\). It follows that
  \(\SQFS \simeq T/J\).
\end{proof}

\end{section}

%\begin{section}{Factorisations and divisibility}
%  \input{factor}
%\end{section}

%\begin{section}{The truncations $\UNI_n$ --- definition and basic properties}
%  \input{truncations}
%\end{section}

%\begin{section}{Presentations of $\UNI_n$ --- defining ideals and the
%    associated simplical complex}
%  \input{presentations}
%\end{section}

%\begin{section}{Homological properties of $\UNI_n$}
%  \input{homological}
%\end{section}

\bibliographystyle{plain}
\bibliography{journals,articles,snellman}
\end{document}